\documentclass[12pt]{article}
\usepackage{amssymb,amsmath,amsopn,amsfonts,mathrsfs,}
\usepackage[dvips]{color}
\usepackage{colordvi,multicol}
\usepackage{epsf}
\usepackage{hyperref}
\usepackage{bbm}
\usepackage[titletoc]{appendix}

\usepackage{xcolor}
\usepackage[T1]{fontenc} % for correct accents in 'Lubomir'

\newcommand{\eps}{\varepsilon}
\newcommand{\ngamma}{{n_\Gamma}}
\newfam\msbfam
\font\tenmsb=msbm10 \textfont\msbfam=\tenmsb \font\sevenmsb=msbm7
\scriptfont\msbfam=\sevenmsb \font\fivemsb=msbm5
\scriptscriptfont\msbfam=\fivemsb

\newcommand{\normmm}[1]{{\left\vert\kern-0.25ex\left\vert\kern-0.25ex\left\vert #1 
    \right\vert\kern-0.25ex\right\vert\kern-0.25ex\right\vert}}

\def\proof{\vspace{1mm}\noindent{\it Proof}\quad}

 \newtheorem{lemma}{Lemma}[section]
 \newtheorem{proposition}[lemma]{Proposition}
 \newtheorem{theorem}[lemma]{Theorem}
 \newtheorem{corollary}[lemma]{Corollary}
 \newtheorem{remark}[lemma]{Remark}

\numberwithin{equation}{section}

  \makeatletter

\newcommand{\Rmnum}[1]{\expandafter\@slowromancap\romannumeral #1@}
\makeatother

\begin{document}
\bigbreak
\title {{\bf Sharp interface limit of stochastic Cahn-Hilliard equation with singular noise
\thanks{Research supported in part  by NSFC (No.11671035). Financial support by the DFG through the CRC 1283 "Taming uncertainty and profiting from randomness and low regularity in analysis, stochastics and their applications" is acknowledged.}\\}}
\author{{\bf \v{L}ubom\'ir Ba\v{n}as}$^{\mbox{c}}$, {\bf Huanyu Yang}$^{\mbox{a,c,d}}$,{\bf Rongchan Zhu}$^{\mbox{b,c}}$,
\date {}
\thanks{E-mail address:  banas@math.uni-bielefeld.de (\v{L}. Ba\v{n}as), hyang@math.uni-bielefeld.de(H. Y. Yang), zhurongchan@126.com(R. C. Zhu)}\\ \\
\small $^{\mbox{a}}$School of Mathematical Science, University of Chinese Academy of Sciences, Beijing 100049, China\\
\small $^{\mbox{b}}$Department of Mathematics, Beijing Institute of Technology, Beijing 100081, China\\
\small $^{\mbox{c}}$ Department of Mathematics, University of Bielefeld, D-33615 Bielefeld, Germany\\
\small $^{\mbox{d}}$ Academy of Mathematics and Systems Science, Chinese Academy of Sciences, Beijing 100190, China}

\maketitle
\begin{abstract}
We study the sharp interface limit of the two dimensional stochastic Cahn-Hilliard equation driven by 
two types of singular noise: a space-time white noise and a space-time singular divergence-type noise.
We show that with appropriate scaling of the noise the solutions of the stochastic problems
converge to the solutions of the determinisitic Mullins-Sekerka/Hele-Shaw problem.
\end{abstract}

\noindent{{\bf Keywords}: stochastic Cahn-Hilliard equation, singular noise, sharp interface limit, Mullins-Sekerka/Hele-Shaw problem.

\section{Introduction}\label{s1}

We consider the stochastic Cahn-Hilliard equation with additive noise
\begin{equation}\label{1.1a}
   \begin{aligned}
   \mathrm{d} u^\varepsilon&=\Delta\Big(-\varepsilon\Delta u^\varepsilon+\frac{1}{\varepsilon}f(u^\varepsilon) \Big)\mathrm{d}t
                                     +\varepsilon^{\sigma}\mathrm{d}\mathcal{W}  \quad &&\mathrm{in}\,\,\mathcal{D}_T:=(0,T)\times\mathcal{D}\,, \\
      \partial_n u^\eps & = \partial_n \Delta u^\varepsilon =0\; &&\text{on} \;(0,T)\times\partial \mathcal{D}\,,\\
   u^\varepsilon(0)&= {u^\eps_0},\\
   \end{aligned}
\end{equation}
where $\mathcal{D}=[0,1]^2$, $n$ is the outward unit normal to $\partial\mathcal{D}$, 
$\sigma>0$ is a constant, $\eps>0$ is a small parameter and $\mathcal{W}$ is a singular space-time noise which will be specified later on.
The nonlinearity in (\ref{1.1a}) is taken as $f(u)=F'(u)$ where $F(u)=\frac{1}{4}(u^2-1)^2$ is the double-well potential.

By introducing an additional variable, a chemical potential $w^\eps$, it is possible to reformulate (\ref{1.1a}) as
\begin{equation}\label{1.1}
   \begin{aligned}
   \mathrm{d} u^\varepsilon&=\Delta w^\eps\mathrm{d}t
                                     +\varepsilon^{\sigma}\mathrm{d}\mathcal{W}  \quad &&\mathrm{in}\,\,\mathcal{D}_T\,, \\
   w^\eps &=-\varepsilon\Delta u^\varepsilon+\frac{1}{\varepsilon}f(u^\varepsilon) \quad &&\mathrm{in}\,\,\mathcal{D}_T\,, \\
      \partial_n u^\eps & = \partial_n w^\varepsilon =0\; &&\text{on} \;(0,T)\times\partial \mathcal{D}\,,\\
   u^\varepsilon(0)&= {u^\eps_0} \quad && \text{in}\; \mathcal{D}.\\
   \end{aligned}
\end{equation}

The deterministic Cahn-Hilliard equation (i.e. (\ref{1.1}) with $\mathcal{W} \equiv  0$) reads as
%\begin{equation}\label{1.2}
%\partial_t u_{\tt D}^\eps=\Delta(-\varepsilon\Delta u_{\tt D}^\eps+\frac{1}{\varepsilon}f(u_{\tt D}^\eps))\, ,
%\end{equation}
\begin{equation}\label{1.4}
%\left\{
   \begin{aligned}
   \partial_t u_{\tt D}^{\varepsilon}&=\Delta w_{\tt D}^\eps   \qquad&&\text{in}\;\mathcal{D}_T,\\
   w_{\tt D}^\eps&=-\varepsilon\Delta u_{\tt D}^{\varepsilon}+\frac{1}{\varepsilon}f(u_{\tt D}^\varepsilon)\; &&\text{in}\;\mathcal{D}_T.\\
   \end{aligned}
%   \right.
\end{equation}
%The Cahn-Hilliard equation is a phenomenological model for the evolution of microstructure in melted binary alloys.
The Cahn-Hilliard equation is a model for the non-equilibrium dynamics of metastable states in phase transitions
\cite{COOK:1970jr,HOHENBERG:1977vo,Langer:1971hu}.
The parameter $\varepsilon$ in (\ref{1.1}) represents an "interaction length", 
which is typically very small, and $u^\eps$ is an order parameter (scaled concentration) 
which assumes the values $u^\eps = 1$ and $u^\eps = - 1$, respectively, in the regions occupied by the pure phases.
The phase separation consists of two stages a so-called spinodal decomposition which is followed by a coarsening process. 
Starting from a fully mixed state, e.g., a random perturbation around the initial mass,
the system undergoes a short phase, so-called spinodal decomposition, during which the initial phases are formed.
The solution quickly approaches the respective values $1$, $-1$ 
in the regions %$\mathcal{D}^+, \, \mathcal{D}^- \subset \mathcal{D}$ 
occupied by the pure phases.
The pure phases are separated by a thin region with a width proportional to $\eps$,
%$\mathcal{D}_\Gamma:=\mathcal{D}\setminus(\mathcal{D}^+\cup\mathcal{D}^-)$ 
so-called diffuse interface. 
Once the diffuse interface is fully formed, the evolution enters a second stage, so-called coarsening phase, during which 
the originally fine-grained structure coarsens, 
the geometric structure of the phase regions gradually becomes simpler and
eventually tends to regions of minimum surface area with preserved volume. 
%In the Cahn-Hilliard equation (\ref{1.1}), 
%this phenomenon corresponds to the behavior of the solution in which the regions $\mathcal{D}^+$, $\mathcal{D}^-$  grow
%and the interfacial region $\mathcal{D}_\Gamma$ eventually tends to a surface having
%minimum surface area (whereas its enclosed region has a fixed volume, since the mass is preserved).
%For more details and discussion we refer to \cite{ALIKAKOS:1994vc} and reference therein.\footnote{will reformulate}

A rigorous sharp interface limit of the deterministic Cahn-Hilliard equation (\ref{1.4}) has been obtained
in \cite{ALIKAKOS:1994vc}
Under the assumption that the interfaces have been formed, i.e., that there exists a smooth closed curve $\Gamma_{00}\subset\subset\mathcal{D}$ 
such that $u^{\varepsilon}(0)\approx -1$ in $\mathcal{D}^-$, in the region enclosed by $\Gamma_{00}$, and $u^{\varepsilon}(0)\approx 1$ in $\mathcal{D}^+:=\mathcal{D}\setminus(\Gamma_{00}\cup\mathcal{D}^-)$,
it is shown in \cite{ALIKAKOS:1994vc} that $v := \lim_{\varepsilon\to 0} w_{\tt D}^\eps$ along with 
$\Gamma_t := \lim_{\varepsilon\to 0} \Gamma_t^{\eps} := \bigl\{ x \in {\mathcal D}:\, u_{\tt D}^{\eps}(t,x) = 0\bigr\}$, $t \in (0, T)$ 
satisfy the deterministic Mullins-Sekerka/Hele-Shaw problem:
\begin{equation}\label{3.2}
\left\{
   \begin{aligned}
   \Delta v&=0 \;\text{in}\;\;\mathcal{D}\setminus\Gamma_t,\;\; t>0,\\
   \partial_nv&=0\;\;\text{on}\;\;\partial\mathcal{D},\\
   v&=\lambda H\;\text{on}\;\;\Gamma_t,\\
   \mathcal{V}&=\frac{1}{2}(\partial_\ngamma v^+-\partial_\ngamma v^-)\;\;\text{on}\;\;\Gamma_t,\\
   \Gamma_0&=\Gamma_{00},\\
   \end{aligned}
   \right.
\end{equation}
where $\lambda = \frac{1}{\sqrt{2}}\int_{-1}^1F(s)\mathrm{d}s$, $H$ is the mean curvature of $\Gamma_t$ with the sign convention that convex hypersurfaces have positive mean curvature, $ \mathcal{V}$ is the normal velocity of the interface with the sign convention that the normal velocity of expanding hypersurfaces is positive, $\ngamma:[0,T]\rightarrow \mathbb{S}^2$ 
is the normal to $\Gamma_t$ and $v^+$, $v^-$ are respectively the restriction of $v$ on {$\mathcal{D}^{\pm}_t$ (the exterior/interior of $\Gamma_t$ in $\mathcal{D}$).}

The sharp interface limit of stochastic Cahn-Hilliard equation with trace-class noise has been studied in \cite{Antonopoulou:2018gh}.
There, the authors show that for sufficiently large $\sigma$ the sharp interface limit of equation (\ref{1.1}) 
satisfies the deterministic Mullins-Sekerka/Hele-Shaw problem (\ref{3.2}). 
In the recent paper \cite{lubo2019} the authors show convergence of structure preserving
numerical approximation of the stochastic Cahn-Hilliard to the deterministic Mullins-Sekerka/Hele-Shaw problem.
In addition \cite{lubo2019} obtains a uniform convergence result which implies convergence of the zero-level set of the numerical solution
to the free boundary $\Gamma:=\cup_{0\leq t\leq T}(\Gamma_t\times\{t\})$ of (\ref{3.2}).
The case of $\sigma\leq 1$ remains an open problem.

In this paper we study the sharp interface limit of stochastic Cahn-Hilliard equation driven by singular noise. 
We consider the Cahn-Hilliard-Cook model, proposed by Cook, cf. \cite{COOK:1970jr} and \cite{HOHENBERG:1977vo}), 
which incorporates thermal fluctuations in the form of an additive noise $\mathcal{W}$ in (\ref{1.1}).
We choose the noise as $\mathcal{W}=W_1$ or 
$\mathcal{W}=\nabla\cdot W_2$, 
where $W_1$ is mass-conserved $L^2(\mathcal{D},\mathbb{R})$-cylindrical Wiener process and $W_2$ is an $L^2(\mathcal{D},\mathbb{R}^2)$-cylindrical Wiener process.
We note that in the latter case  the equation is known as the time-dependent Ginzburg-Landau (TDGL) equation and is
also related to the stochastic quantization for $(\phi)_2^4$-quantum field. 
For the existence and uniqueness results for these two kinds of equations we refer to \cite{DaPrato:1996kk, Rockner:2018ta} and the reference therein.

To analyze the sharp interface limit of the solution $(u^\varepsilon,w^\eps)$ to equation (\ref{1.1}) for the case of the space-time white noise $\mathcal{W}=W_1$,
we adapt the approach of \cite{Antonopoulou:2018gh}.
We estimate the difference of $(u^\varepsilon,w^\eps)$ to an approximate solution $(u^\varepsilon_A,v^A_\varepsilon)$ 
which is constructed by the matched asymptotic expansion method such that the interface $\Gamma_t$ is the zero level set of $u^\varepsilon_A(t)$, cf. \cite{ALIKAKOS:1994vc}.
The approximation satisfies a perturbed equation 
\begin{equation}\label{3.3}
\left\{
\begin{aligned}
\partial_t u^\varepsilon_A&=\Delta w^\eps_A\; &\text{in}\;\mathcal{D}_T,\\
w^\eps_A&=\frac{1}{\varepsilon}f(u^\varepsilon_A)-\varepsilon \Delta u^\varepsilon_A+r^\varepsilon_A\; &\text{in}\;\mathcal{D}_T,\\
\end{aligned}
\right.
\end{equation}
along with the boundary conditions
$$
\frac{\partial u^\varepsilon_A}{\partial n}=\frac{\partial\Delta u^\varepsilon_A}{\partial n}=0\;\text{on}\;\partial \mathcal{D}.
$$
We show that for $\eps \rightarrow 0$ the differences $u^{\varepsilon}-u^\varepsilon_A$, $w^\eps-w^\eps_A$ converge to $0$ for 
$\sigma>\frac{107}{12}$.
As a consequence the sharp interface limit of the equation (\ref{1.1}) satisfies the deterministic Hele-Shaw model (\ref{3.2}) for $\sigma>\frac{107}{12}$. 
We note that the low regularity of the considered noise, prohibits the direct application of It\^o's formula to 
estimate $u^{\varepsilon}-u^\varepsilon_A$. 
Hence the arguments of \cite{Antonopoulou:2018gh} are not directly transferable to our case. 
Instead, we make use of the idea of Da~Prato-Debussche \cite{Debussche:2003du}:
after introducing a variable $Z^{\varepsilon}=\varepsilon^\sigma\int_0^te^{-\varepsilon(t-s)\Delta^2}d\mathcal{W}_s$ 
we study the translated difference $Y^{\varepsilon}=u^{\varepsilon}-u^\varepsilon_A-Z^{\varepsilon}$
which enjoys better regularity properties.
By combining the estimates for $Y^{\varepsilon}$ and $Z^{\varepsilon}$ we we bound the error $u^{\varepsilon}-u^\varepsilon_A$ 
and obtain the sharp-interface limit.
}

For the case singular divergence-type noise $\mathcal{W}=\nabla\cdot W_2$ the equation (\ref{1.1}) is ill-posed in the classical sense, since the solution is not a function but a distribution. 
Hence, it does not make sense to consider the sharp interface limit of (\ref{1.1}) directly. 
Instead we follow the renormalization approach:
we employ a suitable regularization $W^h:=W_2*\rho_h$ of the the space-time white noise $W_2$ and consider the regularized equation:
\begin{equation}\label{1.11}
du^{\varepsilon,h}=\Delta\left(-\varepsilon\Delta u^{\varepsilon,h}+\frac{1}{\varepsilon}\left(f(u^{\varepsilon,h})-3c_{h,t}^{\varepsilon}u^{\varepsilon,h}\right)\right)dt+\varepsilon^{\sigma}\nabla\cdot dW^h_t, 
\end{equation}
where $3c_{h,t}^{\varepsilon}u^{\varepsilon,h}$ is a renormalization term (see Appendix \ref{a.5}) which ensures that $u^{\varepsilon,h}$ converges to $u^{\varepsilon}$ for $h\to 0$,
{ where $u^{\varepsilon}$ is the unique solution of the renormalized version of (\ref{1.1}), see (\ref{5.1}).}
The analysis in the case of the divergence-type noise is complicated by the fact that for fixed $\varepsilon>0$ the renormalization constant in (\ref{1.11}) diverges,
i.e., that $c_{h,t}^{\varepsilon}\to \infty$ as $h\to 0$.
%Consequently it is not obvious how to control the term  $c_{h,t}^{\varepsilon}u^{\varepsilon,h}$.
By choosing $\varepsilon\lesssim h^\iota$ for some $\iota>0$ and $h$ goes to $0$ (see Theorem \ref{t5.4}) the constant 
$c_{h,t}^{\varepsilon}$ becomes small as $h\to 0$ which enables us to control the term $c_{h,t}^{\varepsilon}u^{\varepsilon,h}$. 
The remaining steps in the analysis of (\ref{1.11}) are analogical to the first case:
we obtain that the sharp interface limit of (\ref{1.11}) for $\sigma>\frac{26}{3}$ is the deterministic Mullins-Sekerka problem (\ref{3.2}).

The paper is organized as follows. 
{In \ref{sec_overview} we give an overview of existing results on sharp interface limits for related problems.}
In Section \ref{s2}, we introduce the notation and state preliminary results. 
The sharp interface limit for the space-time white noise is stated in Section \ref{s3} and we prove it in Section \ref{s4}. 
In Section \ref{s5} we use a similar argument as we used in Section \ref{s4} to prove the results for divergence-type noise.

%Background and solution research on deterministic and two types of stochastic equations.
%History of research on Hele-Shaw model.
\section{Overview of existing results}\label{sec_overview}

 For stochastic Cahn-Hilliard eqaution, the authors in \cite{Antonopoulou:2018gh} prove that for large $\sigma$ the sharp interface limit of equation (\ref{1.1}) also satifies the deterministic Hele-Shaw model if $\dot{\mathcal{W}}$ is a trace-class noise. For $\sigma=1$, the sharp interface limit is also conjectured to satisfy the following stochastic Hele-Shaw model:
\begin{equation}
\left\{
   \begin{aligned}
   \Delta v&=0 \;\text{in}\;\mathcal{D}\setminus\Gamma_t,\; t>0,\\
   \partial_nv&=0\;\text{on}\;\partial\mathcal{D},\\
   v&=\lambda H+\mathcal{W}\;\text{on}\;\Gamma_t,\\
    \mathcal{V}&=\frac{1}{2}(\partial_nv^+-\partial_nv^-)\;\text{on}\;\Gamma_t,\\
   \Gamma_0&=\Gamma_{00},\\
   \end{aligned}
   \right.
\end{equation}
In \cite{Antonopoulou:2013tb}, the authors prove that the sharp interface limit of generalized Cahn-Hilliard equation: $\partial_t u=\Delta(-\varepsilon\Delta u+\frac{1}{\varepsilon}f(u)-G_2^{\varepsilon})+G_1^{\varepsilon}$ satisfies the following Hele-Shaw model:
\begin{equation}
\left\{
   \begin{aligned}
   \Delta v&=-\lim_{\varepsilon\to 0}G_1^{\varepsilon} \;\text{in}\;\mathcal{D}\setminus\Gamma_t,\; t>0,\\
   \partial_nv&=0\;\text{on}\;\partial\mathcal{D},\\
   v&=\lambda H-\lim_{\varepsilon\to 0}G_2^{\varepsilon}\;\text{on}\;\Gamma_t,\\
    \mathcal{V}&=\frac{1}{2}(\partial_nv^+-\partial_nv^-)\;\text{on}\;\Gamma_t,\\
   \Gamma_0&=\Gamma_{00},\\
   \end{aligned}
   \right.
\end{equation}
Since they require some regularity conditions for $G_1^\varepsilon$, $G_2^\varepsilon$ w.r.t time, which are not satisfied by Bronwnian motions, it is not clear how to obtain the stochastic Hele-Shaw model rigorously. Until now, the rigorous complete description of the motion of interfaces in dimensions two and three in stochastic case stands for many years as a wide open problem.

Another simpler model is the following Allen-Cahn equation
\begin{equation}\label{1.3}
\partial_t u=\Delta u-\frac{1}{\varepsilon^2}f(u).
\end{equation}
It is well-known that the movement of interface is characterized by mean curvature flow (see e.g.\cite{Evans:1992cc,ILMANEN:1993uu,deMottoni:1995ff}). Unlike the solution to the Allen-Cahn equation,
the solution to the Cahn-Hilliard equation (\ref{1.4}) does not approach $\pm 1$ away from the interface exponentially fast. The direct application of the method of asymptotic matching in \cite{DEMOTTONI:1995una} does not lead to the desired approximation solutions. 
In stochastic case which is also called Model A of \cite{HOHENBERG:1977vo}),  the authors in \cite{FUNAKI:1999ho} and \cite{Weber:2010dl} consider the following stochastic Allen-Cahn equation
\begin{equation}
\partial_tu=\Delta u-\frac{1}{\varepsilon^2}f(u)+\frac{1}{\varepsilon}\Xi_t^\varepsilon.
\end{equation}
The noise $\Xi^\varepsilon$ is constant in space and smooth in time. For $\varepsilon\to 0$
the correlation length goes to zero at a precise rate and  $
\int_0^t\Xi^\varepsilon_sds$ converges to a Brownian motion pathwisely. They prove that the dynamics of the phase-separating hyperplane  $\Gamma_t$ appearing in the limit is given by stochastic mean curvature flow (see also in \cite[Chapter 4]{Funaki:2016iz}). For space-time white noise, in \cite{Tsatsoulis:2018ws} the authors prove the "exponential loss of memory property". But for sharp interface limit, there is still no result for space-time white noise.

\section{Notations and preliminaries}\label{s2}

Throughout the paper, we use the notation $a\lesssim b$ if there exists a constant $c>0$ which is independent to $\varepsilon$ and time $T$ 
such that $a\leq cb$. If $c$ is depend on $T$, we use the notation $a\lesssim_T b$.
We write $a\backsimeq b$ if $a\lesssim b$ and $b\lesssim a$.

Let $\mathcal{D}:=(0,1)^2$, $\mathcal{D}_T:=(0,T)\times\mathcal{D}$. In this paper, we always use $\langle\cdot,\cdot\rangle$ to denote the $L^2(\mathcal{D})$-inner product. For any $E\subset\mathcal{D}$, we denote by $\chi_E$ the characteristic function of $E$, i.e.
\begin{equation*}
\chi_E(x)=
\left\{
\begin{aligned}
1\;\;\text{if}\;\;x\in E,\\
0\;\;\text{if}\;\;x\not\in E.
\end{aligned}
\right.
\end{equation*}

We consider the Neumann Laplacian operator $\Delta$ on $L^2(\mathcal{D})$ with domain
$$
D(\Delta)=\{u\in H^2(\mathcal{D}):\frac{\partial u}{\partial n}=0\;\text{on}\;\partial \mathcal{D} \}.
$$
The operator $-\Delta$ is self-adjoint positive and has compact resolvent. It possesses a basis of eigenvectors $\{e_k\}_{k\in \mathbb{Z}^2}$ which is orthonormal in $L^2(\mathcal{D})$. In fact for $k=(k_1,k_2)\in\mathbb{Z}^2$, $e_k(x)$ is given by 
\begin{equation}\label{2.3}
 \begin{aligned}
 & e_0(x):=1,e_{(k_1,0)}(x)=\sqrt{2}\cos \pi k_1x_1, e_{(0,k_2)}(x)=\sqrt{2}\cos \pi k_2x_2,,\\ &e_k(x):=2\cos \pi k_1x_1\cdot\cos \pi k_2x_2,k_1k_2\neq 0.
\end{aligned}
\end{equation}
It is associated with the eigenvalues $\{\lambda_k\}$, where $\lambda_k\simeq |k|^2$. 

We also introduce a notation for the average of $g\in L^2(\mathcal{D})$:
  $$m(g):=\langle g,e_0\rangle.$$
  
  For any $\alpha\in \mathbb{R}$, we define $V^{\alpha}$ as the closure of $C^\infty(\mathcal{D})$ under the norm
  $$
\|g\|_{V^\alpha}^2:=m(g)^2+\sum_{k}\lambda_k^{\alpha}\langle g,e_k \rangle^2.  
  $$
It is easy to see that $(V^{\alpha},\|\cdot\|_{V^\alpha})$ is a Hilbert space and $V^{\alpha}\simeq H^{\alpha}$, where $H^\alpha$ is the classical Sobolev space on domain $\mathcal{D}$ which can be defined as the closure of $C^\infty(\mathcal{D})$ under the norm
$$
\|g\|_{H^\alpha}^2=\sum_{k\in\mathbb{Z}^2}(1+\lambda_k)^\alpha\langle g,e_k\rangle^2.
$$ 

In the rest of this paper, we use the notation $H^\alpha$ to represent $V^\alpha$ for simiplicity.

Moreover for any $s,\alpha\in \mathbb{R}$, we can define a bounded operator $(-\Delta)^{s}:H^\alpha\to H^{\alpha-2s}$ by:
$$(-\Delta)^{s}u=\sum_{k\in\mathbb{Z}^2\setminus\{(0,0)\}}\lambda_k^{s}u_ke_k,$$
where $u=\sum_{k}u_ke_k\in H^{\alpha}$.

We also set
  $$H_0^{\alpha}:=\{g \in H^{\alpha}:\langle g,e_0\rangle _{H^{\alpha}}=0\},$$
where $\langle \cdot,\cdot\rangle _{H^{\alpha}}$ denote the inner product in $H^\alpha$.  Moreover we denote $L_0^2:=H_0^0$. 

The analysis of this paper relies heavily on the existence of smooth solution to (\ref{3.2}) which is guaranteed by the next theorem.
\begin{theorem}{\cite[Theorem 1.1]{Xinfu:1996ci}}\label{t2.1}
For any $\Gamma_{00}\in\mathcal{C}^{3+\alpha}$ for some $\alpha \in(0,1)$, there exists a $T>0$, such that (\ref{3.2}) has a unique local solution $\{(v,\Gamma)\}_{t\in[0,T]}$,
where $\Gamma\in C^{\frac{3+\alpha}{3}}([0,T];\mathcal{C}^{3+\alpha})$.
\end{theorem}
Throughout the paper we assume that $\Gamma_{00}$ and $T$ 
satisfy the conditions of Theorem~\ref{t2.1}, i.e., that the Mullins-Sekerka problem (\ref{3.2}) admits a unique classical solution  on $[0,T]$.
Consequently, it is possible to construct an approximate solution that satisfies (\ref{3.3}).
The properties of the solution of (\ref{3.3}) which are summarized in the theorem below {are the consequence of Theorems~2.1~and~4.12 of \cite{ALIKAKOS:1994vc}}.
\begin{theorem}\label{t3.2}
Let $\{(v,\Gamma_t)\}_{t\in[0,T]}$ be the classical smooth solution to (\ref{3.2}).
For any $K>0$ there exists a pair $(u^\varepsilon_A,w^\eps_A)$ of solutions to (\ref{3.3}), such that for a small enough $\varepsilon>0$, $u^\varepsilon_A$ is uniformly bounded and
$$
\|r^\varepsilon_A\|_{C(\mathcal{D}_T)}\lesssim \varepsilon^{K-2},
$$
Moreover
$$
\|w^\eps_A-v\|_{C(\mathcal{D}_T)}\lesssim\varepsilon,
$$
where $v$ is the solution to (\ref{3.2}) below.

{ Finally for any $x\in \mathcal{D}\setminus\Gamma_t$ such that $d(x,\Gamma_t)>C \varepsilon^\alpha$, where $0<\alpha<1$ and $d(x,\Gamma_t)$
is the distance of $x$ to $\Gamma_t$ and $C$ is a constant that is independent to $\varepsilon$ ,
it holds respectively that $|u^\varepsilon_A(t,x)-1|\lesssim \varepsilon$ on $\mathcal{D}^+_t$, $|u^\varepsilon_A(t,x)+1|\lesssim \varepsilon$
on $\mathcal{D}^-_t$, where $\mathcal{D}^-_t$ and $\mathcal{D}^+_t$ are respectively the interior and exterior of
$\Gamma_t$ in $\mathcal{D}$. $u^\varepsilon_A$ also satisfies the following thin interface conditions:
$$
m\left(\{(t,x)\in\mathcal{D}_T:|u^\varepsilon_A|<1\}\right)\lesssim\varepsilon.
$$
}
\end{theorem}

\section{The sharp interface limit for space-time white noise}\label{s3} 

 Let $\mathcal{W}=W$ be an $L^2_0(\mathcal{D})$-cylindrical Wiener process on a fixed stochastic basis $\left(\Omega,\mathcal{F},\mathbb{P}\right)$. 
\begin{theorem}\label{t3.1}
(\cite[Theorem 2.1]{DaPrato:1996kk})For $\mathbb{P}-a.s.\;\omega$,
there exists a unique solution $u^{\varepsilon}$ to equation (\ref{1.1}) in $C([0,T];H^{-1})$.
\end{theorem}

We rewrite the equation (\ref{1.1}) as
\begin{equation}\label{3.1}
\left\{
   \begin{aligned}
   \mathrm{d}u^{\varepsilon}&=\Delta w^\eps\mathrm{d}t+\varepsilon^{\sigma}dW \;\text{in}\;\mathcal{D}_T,\\
   w^\eps&=\frac{1}{\varepsilon}f(u^\varepsilon)-\varepsilon\Delta u^{\varepsilon}\;\text{in}\;\mathcal{D}_T.\\
   \end{aligned}
   \right.
\end{equation}
We assume that the interface has been formed initially. That is, there exists a smooth closed curve $\Gamma_{00}\subset\subset\mathcal{D}$ such that $u^{\varepsilon}(0)\approx -1$ in $\mathcal{D}^-$, the region enclosed by $\Gamma_{00}$, and $u^{\varepsilon}(0)\approx 1$ in $\mathcal{D}^+:=\mathcal{D}\setminus(\Gamma_{00}\cup\mathcal{D}^-)$.

Our main theorem will show that as $\varepsilon\to 0$, $w^\eps$ tends to $v$, which, together with a free boundary $\Gamma\equiv\cup_{0\leq t\leq T}(\Gamma_t\times\{t\})$, satisfies the deterministic Hele-Shaw problem (\ref{3.2}).

We present now the following spectral estimate which is useful in our proof.
\begin{proposition}\label{p3.3}(\cite[Proposition 3.1]{ALIKAKOS:1994vc})
Let $u^\varepsilon_A$ be the approximation given in Theorem \ref{t3.2}. Then for all $w\in H^1$ satisfying Neumann boundary conditions such that $\int_{\mathcal{D}}w=0$, the following estimate is valid
$$
\varepsilon\|w\|_{H^1}^2+\frac{1}{\varepsilon}\int f'(u^\varepsilon_A)w^2\geq -C_0\|w\|_{H^{-1}}^2.
$$
\end{proposition}

We consider the residual
\begin{equation}\label{3.4}
R^\varepsilon:=u^{\varepsilon}-u^\varepsilon_A,
\end{equation}
where $u^{\varepsilon}$ is the unique solution to (\ref{3.1}). We show bounds for this error $R^\varepsilon$ in our main theorem below.
\begin{theorem}\label{t3.4}(\textbf{Main Theorem})
Let $u^\varepsilon_A$ be defined in Theorem \ref{t3.2} with large enough $K$ and let $u^{\varepsilon}$ be the unique solution to (\ref{1.1}) with initial value $u^{\varepsilon}(0)=u^\varepsilon_A(0)$. For any $\sigma^*>\delta>0$,
\begin{equation*}
\left\{
\begin{aligned}
\gamma&>13,\\
\sigma^\ast&>\frac{1}{3}\gamma+\frac{13}{3}+2\delta,
\end{aligned}\right.
\end{equation*}
where $\sigma^\ast=\sigma-\frac{1}{4}$ is introduced in Lemma \ref{l4.1},
there exist a generic constant $C>0$ and a constant $C_\delta>0$ for all $\delta>0$ such that the following estimates hold
\begin{equation*}
\begin{aligned}
\mathbb{P}\left[\|R^\varepsilon\|_{L^3(\mathcal{D}_T)}\leq C\varepsilon^{\frac{\gamma}{3}}\right]\geq 1-C_\delta \varepsilon^\delta,\\
\mathbb{P}\left[\|R^\varepsilon\|_{L^{\infty}(0,T;H^{-1})}^2\leq C\left(\varepsilon^{\gamma-1}+\varepsilon^{\sigma^\ast-1-2\delta+\frac{\gamma}{3}}\right)\right]\geq 1-C_\delta\varepsilon^\delta,\\
\mathbb{P}\left[\|w^\eps-w^\eps_A\|_{L^{1}(0,T;H^{-2})}^2\leq C\varepsilon^{\frac{\gamma}{3}-1}\right]\geq 1-C_\delta\varepsilon^\delta.
\end{aligned}
\end{equation*}
\end{theorem}

\begin{remark}
Since $\delta$ can be arbitrarily small, the best choice is $\sigma>\frac{107}{12}$.
\end{remark}

\begin{corollary}\label{c3.5}
There exists a subsequence $\{\varepsilon_k\}_{k=1}^\infty$ such that for $\mathbb{P}-a.s.\;\omega\in\Omega$
$$
\lim_{k\to\infty}u^{\varepsilon_k}={ 1-2\chi_{\mathcal{D}^-_t}}\;\;\text{in}\;\;L^3({\mathcal{D}_{T}}),
$$
where $\mathcal{D}^-_t$ is the interior of $\Gamma_t$ in $\mathcal{D}$.
\end{corollary}
\proof
%The local uniqueness of (\ref{3.2}) can be obtained directly by \cite[Theorem 1.1]{Xinfu:1996ci}. 
{We note that by Theorem \ref{t3.4} the problem (\ref{3.2}) has a unique strong solution on $[0,T]$.}
Hence, by the construction of $u^\varepsilon_A$, see \cite{ALIKAKOS:1994vc}, it holds uniformly $t\in[0,T]$ that
$$
\lim_{\varepsilon\to 0}u^\varepsilon_A=1-2\chi_{\mathcal{D}^-_t}\;\;\text{uniformly on compact subsets}.
$$

For any $\eta>0$, choosing $\varepsilon$ small enough such that $C\varepsilon^{\frac{\gamma}{3}}<\eta$, then we have
$$
\mathbb{P}\left[\|R^\varepsilon\|_{L^3(\mathcal{D}_T)}>\eta\right]\leq \mathbb{P}\left[\|R^\varepsilon\|_{L^3(\mathcal{D}_T)}> C\varepsilon^{\frac{\gamma}{3}}\right]\leq C_\delta \varepsilon^\delta,
$$
which implies that $\|R^\varepsilon\|_{L^3}$ converge in probability to $0$. Thus there exists a subsequence (still denoted as $\varepsilon$), such that 
$$
\lim_{\varepsilon\to 0}\|R^{\varepsilon}\|_{L^3(\mathcal{D}_{T})}=0\;\;\mathbb{P}-a.s..
$$
Since $R^\varepsilon=u^\varepsilon-u^\varepsilon_A$, we obtain the assertion.

$\hfill\Box$
\vskip.10in
\subsection{The proof of the Main Theorem}\label{s4}
\subsubsection{The decomposition of the equation for the error}
On Combining (\ref{3.1}), (\ref{3.3}) and noting (\ref{3.4}) we obtain
\begin{equation}\label{4.1}
\left\{
\begin{aligned}
dR^\varepsilon&=-\varepsilon \Delta^2R^\varepsilon dt+\frac{1}{\varepsilon}\Delta\left(f(u^\varepsilon_A+R^\varepsilon)-f(u^\varepsilon_A)\right)dt+\Delta r^\varepsilon_Adt+\varepsilon^{\sigma}dW,\\
\frac{\partial R^{\varepsilon}}{\partial n}&=\frac{\partial \Delta R^{\varepsilon}}{\partial n}=0\;\;\text{on}\;\;\partial\mathcal{D}.
\end{aligned}\right.
\end{equation}
Let $Z^{\varepsilon}_t:=\varepsilon^{\sigma}\int_0^te^{-(t-s)\varepsilon\Delta^2}dW_s$, which is the mild solution to the linear equation:
\begin{equation}\label{4.2}
\left\{
\begin{aligned}
dZ^{\varepsilon}&=-\varepsilon\Delta^2Z^{\varepsilon}dt+\varepsilon^{\sigma}dW,\\
\frac{\partial Z^{\varepsilon}}{\partial n}&=\frac{\partial \Delta Z^{\varepsilon}}{\partial n}=0\;\;\text{on}\;\;\partial\mathcal{D}.
\end{aligned}\right.
\end{equation}
Then $Y^{\varepsilon}:=R^\varepsilon-Z^{\varepsilon}$ satisfies:
\begin{equation}\label{4.3}
\left\{
\begin{aligned}
dY^{\varepsilon}&=-\varepsilon\Delta^2Y^{\varepsilon}dt+\frac{1}{\varepsilon}\Delta\left(f'(u^\varepsilon_A)(Y^{\varepsilon}+Z^{\varepsilon})+\mathcal{N}(u^\varepsilon_A,Y^{\varepsilon}+Z^{\varepsilon})\right)dt+\Delta r^\varepsilon_Adt,\\
\frac{\partial Y^{\varepsilon}}{\partial n}&=\frac{\partial \Delta Y^{\varepsilon}}{\partial n}=0\;\;\text{on}\;\;\partial\mathcal{D}.
\end{aligned}\right.
\end{equation}
where $\mathcal{N}(u,v):=f(u+v)-f(u)-f'(u)v$.

Moreover, we define a stopping time $T_{\varepsilon}$ by:
\begin{equation}\label{4.4}
T_{\varepsilon}:=T\wedge \inf\{t>0:\int_0^t\|Y^{\varepsilon}_s\|_{L^3}^3ds>\varepsilon^{\gamma}\},
\end{equation}
for some $\gamma>1$.

\subsubsection{Estimate for $Z^{\varepsilon}$}

\begin{lemma}\label{l4.1}
For any $\delta>0$, there exists a constant $C_\delta>0$, such that
$$\mathbb{P}[\Omega_{\delta}]>1-C_\delta\varepsilon^\delta,$$
where $C_1>0$ is a universal constant, $\Omega_{\delta}:=\{\|Z\|_{C(\mathcal{D}_T)}\leq C_1\varepsilon^{\sigma^\ast-2\delta}\}$, and $\sigma^\ast:=\sigma-\frac{1}{4}$.
\end{lemma}
\proof By the factorization method in \cite{DaPrato:2004fs} we have that for $\kappa\in(0,1)$
$$Z^{\varepsilon}(t)=\varepsilon^{\sigma}\frac{\sin(\pi\kappa)}{\pi}\int_0^t(t-s)^{\kappa-1}\langle M(\varepsilon(t-s),x,\cdot),U(s)\rangle ds,$$
where $M(\varepsilon t,x,y)$ is the kernel of the semigroup $\{e^{-\varepsilon t\Delta^2}\}$
and
$$U^{\varepsilon}(s,\cdot)=\int_0^s(s-r)^{-\kappa} e^{-\varepsilon(s-r) \Delta^2}dW_r.$$
Similarly to the proof of Lemma 2.12 in \cite{DaPrato:2004fs}, we have that
\begin{equation}\label{2.4}
\mathbb{E}\left[\|Z^{\varepsilon}(t)\|_{\mathcal{C}(\mathcal{D}_T)}\right]\lesssim_T\varepsilon^\sigma\mathbb{E}\left[\|U^{\varepsilon}\|_{L^{2p}(\mathcal{D}_T)}\right].
\end{equation}
It suffices to estimate $\mathbb{E}\left[\|U^{\varepsilon}\|_{L^{2p}(\mathcal{D}_T)}\right]$ for $p>\frac{1}{2\kappa}$.

In fact, we have that
\begin{equation}\label{2.5.0}
\begin{aligned}\mathbb{E}\left[\|U^{\varepsilon}(s)\|_{L^{2p}(\mathcal{D}_T)}^{2p}\right]\lesssim 
&\int_{\mathcal{D}_T}\mathbb{E}\left[\left|\int_0^s(s-r)^{-\kappa} e^{-\varepsilon(s-r)\Delta^2}dW_r\right|^{2p}\right]dsdx\\
\lesssim&\int_{\mathcal{D}_T}\left(\mathbb{E}\left[\left|\int_0^s(s-r)^{-\kappa} e^{-\varepsilon(s-r)\Delta^2} dW_r\right|^{2}\right]\right)^p dsdx.\
\end{aligned}
\end{equation}
Here we used that $U^\varepsilon(x)$ belongs to the first order Wiener-chaos and Gaussian hypercontractivity (cf. \cite[Section 1.4.3]{Nualart:1995er} and \cite{thefreemarkofffie:4RNWhXc2}) in the second inequality.
Moreover, we obtain that
\begin{equation}\label{2.5}
\begin{aligned}
\mathbb{E}\left[\left|\int_0^s(s-r)^{-\kappa} e^{-\varepsilon(s-r)\Delta^2}dW_r\right|^{2}\right]
\lesssim \int_0^s\int_{\mathcal{D}}(s-r)^{-2\kappa}M(\varepsilon(s-r),x,y)^2dyds.
\end{aligned}
\end{equation}
Since $M(t,x,y)$ is the kernel of $e^{-t\Delta^2}$, we have that for any $g\in L^2$
\begin{equation*}
\int_{\mathcal{D}}M(t,x,y)g(y)dy =e^{-t\Delta^2}g(x)\simeq\sum_k\langle g,e^{- t|k|^4}e_k\rangle e_k(x).
\end{equation*}
Hence 
\begin{equation}\label{2.6}
M(t,x,y)\simeq\sum_k e^{-t|k|^4}e_k(x)e_k(y).
\end{equation}
where $e_k$ is defined in (\ref{2.3}). Note that $e_k(x)e_k(y)=\frac{1}{2}(e_k(x-y)+e_k(x+y)$. Thus we obtain 
\begin{equation}\label{2.6.1}
M(t,x,y)\simeq\sum_k e^{-t|k|^4}(e_k(x-y)+e_k(x+y)):=P(t,x-y)+P(t,x+y),
\end{equation}
Then (\ref{2.5}) becomes
\begin{equation}\label{2.7}
\mathbb{E}\left[\left|\int_0^s(s-r)^\kappa e^{-\varepsilon(s-r)\Delta^2} dW_r\right|^{2}\right]
\lesssim \int_0^s\int_{\mathcal{D}}(s-r)^{-2\kappa}\left(P(\varepsilon(s-r),x-y)^2+P(\varepsilon(s-r),x+y)^2\right)dyds,.
\end{equation}
By \cite[p282, (c)]{Stein:1972jf}, we have that
\begin{equation}\label{2.8}
|P(t,x)|\lesssim |x|^{-2}e^{-\frac{t}{|x|^4}}\lesssim t^{-\frac{\eta}{4}}|x|^{-2+\eta},\;\forall \eta\in[0,2].
\end{equation}
Then taking (\ref{2.7}) into (\ref{2.8}), we deduce that 
\begin{equation}\label{2.9}
\begin{aligned}
\mathbb{E}\left[\left|\int_0^s(s-r)^\kappa e^{-\varepsilon(s-r)\Delta^2}dW_r\right|^{2}\right]
&\lesssim \varepsilon^{-\frac{\eta}{2}}\int_0^s\int_{\mathcal{D}}(s-r)^{-2\kappa-\frac{\eta}{2}}\left(|x+y|^{-4+2\eta}+|x-y|^{-4+2\eta}\right)dyds\\
&\lesssim \varepsilon^{-\frac{\eta}{2}}s^{1-2\kappa-\frac{\eta}{2}}|x|^{-2+2\eta}.
\end{aligned}
\end{equation}
Here we require that
$$
1-2\kappa-\frac{\eta}{2}>0,\quad -2+2\eta>0,
$$
that is
\begin{equation}\label{2.10}
1<\eta<2-4\kappa,
\end{equation}
which can be obtained by choosing small enough $\kappa>0$.
Hence by (\ref{2.4}) and (\ref{2.5.0}), we obtain that for any $p\geq 1$
$$\mathbb{E}\left[\|U^{\varepsilon}\|_{L^{2p}(\mathcal{D}_T)}\right]\lesssim \varepsilon^{\sigma-\frac{\eta}{4}},$$
This implies that for any $2>\eta>1$, 
\begin{equation}
\mathbb{E}\left[\|Z^\varepsilon\|_{\mathcal{C}(\mathcal{D}_T)}\right]\lesssim\varepsilon^{\sigma-\frac{\eta}{4}}.
\end{equation}
Hence the statement follows by Cheybeshev's inequatliy.
  
$\hfill\Box$
\vskip.10in

\subsubsection{Local-in-time estimate for $Y^{\varepsilon}$ up to $T_{\varepsilon}$ on the set $\Omega_\delta$} 
In the remainder of the proof we fix $\omega\in \Omega_{\delta}$ where $\Omega_{\delta}$ is defined in Lemma \ref{l4.1}
and work pathwise.
We note that by definition $\|Z^{\varepsilon}(\omega)\|_{C(\mathcal{D}_T)}\lesssim \varepsilon^{\sigma^\ast-2\delta}$ for $\omega\in \Omega_{\delta}$.

By taking inner product with $(-\Delta)^{-1}Y^{\varepsilon}$ in both side of equation (\ref{4.3}) we have that
\begin{equation}\label{4.5}
\frac{1}{2}\frac{d\|Y^{\varepsilon}\|_{H^{-1}}^2}{dt}+\varepsilon\|Y^{\varepsilon}_t\|_{H^1}^2=-\frac{1}{\varepsilon}\langle f'(u^\varepsilon_A)(Y^{\varepsilon}+Z^{\varepsilon})+\mathcal{N}(u^\varepsilon_A,Y^{\varepsilon}+Z^{\varepsilon}),Y^{\varepsilon}\rangle-\langle r^\varepsilon_A,Y^{\varepsilon}\rangle.
\end{equation}
We estimate the right hand side of (\ref{4.5}) separately.
Using Proposition \ref{p3.3} we have that
\begin{equation}\label{4.6}
-\frac{1}{\varepsilon}\langle f'(u^\varepsilon_A)Y^{\varepsilon},Y^{\varepsilon}\rangle 
\leq\varepsilon\|Y^\varepsilon\|_{H^1}^2+C_0\|Y^{\varepsilon}\|_{H^{-1}}^2
\end{equation}
For $-\frac{1}{\varepsilon}\langle f''(u^\varepsilon_A)(Y^{\varepsilon},Z^{\varepsilon}\rangle $ by Theorem \ref{t3.2} we know that $u^\varepsilon_A$ is uniformly bounded in $\mathcal{D}_T$. Thus we have that
\begin{equation}\label{4.7}
\frac{1}{\varepsilon}|\langle f''(u^\varepsilon_A)Y^{\varepsilon},Z^{\varepsilon}\rangle|
\lesssim \frac{1}{\varepsilon}\|Y^{\varepsilon}\|_{L^3}\|Z^{\varepsilon}\|_{L^{\frac{3}{2}}}
\lesssim\varepsilon^{\sigma^\ast-1-2\delta}\|Y^{\varepsilon}\|_{L^3},
\end{equation}
where we used H\"older's inequality in the first inequality and Lemma \ref{l4.1} in the last inequality.

By \cite[Lemma 2.2]{ALIKAKOS:1994vc}, we have that
$
v\mathcal{N}(u,v)\geq-C|v|^3.
$
Then
\begin{equation}\label{4.8}
\begin{aligned}
-\frac{1}{\varepsilon}\langle \mathcal{N}(u^\varepsilon_A,Y^{\varepsilon}+Z^{\varepsilon}),Y^{\varepsilon}\rangle
&= -\frac{1}{\varepsilon}\langle \mathcal{N}(u^\varepsilon_A,Y^{\varepsilon}+Z^{\varepsilon}),Y^{\varepsilon}+Z^{\varepsilon}\rangle+\frac{1}{\varepsilon}\langle \mathcal{N}(u^\varepsilon_A,Y^{\varepsilon}+Z^{\varepsilon}),Z^{\varepsilon}\rangle\\
&\lesssim\frac{1}{\varepsilon}\|Y^{\varepsilon}+Z^{\varepsilon}\|_{L^3}^3+\frac{1}{\varepsilon}|\langle \mathcal{N}(u^\varepsilon_A,Y^{\varepsilon}+Z^{\varepsilon}),Z^{\varepsilon}\rangle|\\
&\lesssim \frac{1}{\varepsilon}\|Y\|_{L^3}^3+\varepsilon^{3(\sigma^\ast-2\delta)-1}+\frac{1}{\varepsilon}|\langle \mathcal{N}(u^\varepsilon_A,Y^{\varepsilon}+Z^{\varepsilon}),Z^{\varepsilon}\rangle|,
\end{aligned}
\end{equation}
where we used Lemma \ref{l4.1} in the last inequality.

For $|\langle \mathcal{N}(u^\varepsilon_A,Y^{\varepsilon}+Z^{\varepsilon}),Z^{\varepsilon}\rangle|$, by the Taylor expansion, $\mathcal{N}(u^\varepsilon_A,Y^{\varepsilon}+Z^{\varepsilon})=f''(u^\varepsilon_A+\theta(Y^{\varepsilon}+Z^{\varepsilon}))(Y^{\varepsilon}+Z^{\varepsilon})^2=6(u^\varepsilon_A+\theta(Y^{\varepsilon}+Z^{\varepsilon}))(Y^{\varepsilon}+Z^{\varepsilon})^2$, where $\theta\in (0,1)$. Then we have
\begin{equation}
\begin{aligned}\label{4.9}
|\langle \mathcal{N}(u^\varepsilon_A,Y^{\varepsilon}+Z^{\varepsilon}),Z^{\varepsilon}\rangle|
&\lesssim \varepsilon^{\sigma^\ast-2\delta}\|\mathcal{N}(u^\varepsilon_A,Y^{\varepsilon}+Z^{\varepsilon})\|_{L^1}\\
&\lesssim \varepsilon^{\sigma^\ast-2\delta}(\|Y^{\varepsilon}+Z^{\varepsilon}\|_{L^3}^3+\|Y^{\varepsilon}+Z^{\varepsilon}\|_{L^2}^2)\\
&\lesssim \varepsilon^{3(\sigma^\ast-2\delta)}+\varepsilon^{4(\sigma^\ast-2\delta)}+\varepsilon^{\sigma^\ast-2\delta}\|Y^{\varepsilon}\|_{L^3}^2+\varepsilon^{\sigma^\ast-2\delta}\|Y^{\varepsilon}\|_{L^3}^3,
\end{aligned}
\end{equation}
where we used the uniform boundness of $u^\varepsilon_A$ in the second inequality and Lemma \ref{l4.1} in the first and the last inequality.

For $|\langle r^\varepsilon_A, Y^{\varepsilon}\rangle|$, by Theorem \ref{t3.2} we have
\begin{equation}\label{4.10}
|\langle r^\varepsilon_A, Y^{\varepsilon}\rangle|\lesssim\varepsilon^{K-2}\|Y^{\varepsilon}\|_{L^1}\lesssim \varepsilon^{K-2}\|Y^{\varepsilon}\|_{L^3}.
\end{equation}

Let $\sigma^\ast>\delta$, $\varepsilon<1$, $\delta$ be small enough and $K$ large enough. Collecting (\ref{4.5})-(\ref{4.10}) together,  by using H\"older's inequality we have

$$
\frac{d\|Y^{\varepsilon}(t)\|_{H^{-1}}^2}{dt}\lesssim \|Y^{\varepsilon}\|_{H^{-1}}^2+
\frac{1}{\varepsilon}\|Y^{\varepsilon}\|_{L^3}^3+\varepsilon^{\sigma^\ast-1-2\delta}(\|Y^{\varepsilon}\|_{L^3}+\|Y^{\varepsilon}\|_{L^3}^2+\|Y^{\varepsilon}\|_{L^3}^3)+\varepsilon^{3(\sigma^\ast-2\delta)-1}.
$$
Then for any $t\leq T_{\varepsilon}$ we have
\begin{equation}\label{4.11}
\begin{aligned}
\|Y^{\varepsilon}(t)\|_{H^{-1}}^2&\lesssim \int_0^te^{t-s}\left(\frac{1}{\varepsilon}\|Y^{\varepsilon}\|_{L^3}^3+\varepsilon^{\sigma^\ast-1-2\delta}\|Y^{\varepsilon}\|_{L^3}+\varepsilon^{3(\sigma^\ast-2\delta-)-1}\right)ds\\
&\lesssim_T \frac{1}{\varepsilon}\int_0^t\|Y^{\varepsilon}\|_{L^3}^3d\tau+ \varepsilon^{\sigma^\ast-1-2\delta}\left(\int_0^t\|Y^{\varepsilon}\|_{L^3}^3d\tau\right)^{\frac{1}{3}}+\varepsilon^{3(\sigma^\ast-2\delta)-1}\\
&\lesssim \varepsilon^{\gamma-1}+\varepsilon^{\sigma^\ast-1-2\delta+\frac{\gamma}{3}}+\varepsilon^{3(\sigma^\ast-2\delta)-1}.
\end{aligned}
\end{equation}

To estimate $L^2(0,T_{\varepsilon};H^1)$ norm of $Y^{\varepsilon}$, we use the estimate presented in \cite[p.171]{ALIKAKOS:1994vc}
$$
-\frac{1}{\varepsilon}\int_0^t\int_{\mathcal{D}}f''(u^\varepsilon_A)g^2dxds\lesssim \varepsilon^{-\frac{2}{3}}(\int_0^t\|g\|_{L^3}^3ds)^{\frac{2}{3}},\;\;\forall g\in L^3.
$$
Then
\begin{equation}\label{4.12}
-\frac{1}{\varepsilon}\int_0^t\langle f''(u^\varepsilon_A)Y^{\varepsilon},Y^{\varepsilon}\rangle ds
\leq\varepsilon^{-\frac{2}{3}}(\int_0^t\|Y^{\varepsilon}\|_{L^3}^3ds)^{\frac{2}{3}}\lesssim \varepsilon^{\frac{2}{3}(\gamma-1)}.
\end{equation}

Combining (\ref{4.5}), (\ref{4.7})-(\ref{4.10}) and (\ref{4.12}) we have for any $t\leq T_{\varepsilon}$
\begin{equation}\label{4.13}
\int_0^t\|Y^{\varepsilon}\|_{H^1}^2ds\lesssim \varepsilon^{\frac{2}{3}\gamma-\frac{5}{3}}+\varepsilon^{\sigma^\ast-2-2\delta+\frac{\gamma}{3}}+\varepsilon^{3(\sigma^\ast-2\delta)-2}+\varepsilon^{\gamma-2}.
\end{equation}

\subsubsection{Final step: Globalization $T_\varepsilon\equiv T$}
Let $$\gamma_1:=(\gamma-1)\wedge (3(\sigma^\ast-2\delta)-1)\wedge (\sigma^\ast-1-2\delta+\frac{\gamma}{3}),$$ 
$$\gamma_2:=(\frac{2}{3}\gamma-\frac{5}{3})\wedge(3(\sigma^\ast-2\delta)-2)\wedge(\sigma^\ast-2-2\delta+\frac{\gamma}{3})\wedge(\gamma-2)=(\frac{2}{3}\gamma-\frac{5}{3})\wedge(\gamma_1-1),$$ 
then we have for any $t\leq T_{\varepsilon}$
\begin{equation}\label{4.14}
\sup_{s\in[0,t]}\|Y^{\varepsilon}\|_{H^{-1}}^2\lesssim \varepsilon^{\gamma_1},\quad \int_0^t\|Y^{\varepsilon}\|_{H^1}^2ds\lesssim\varepsilon^{\gamma_2}.
\end{equation}

We use the Sobolev's embedding of $H^{\beta}$ into $L^p$ with $\beta:=2(\frac{1}{2}-\frac{1}{p})=\frac{p-2}{p}$.  Then by the interpolation we have 
$$
\|Y^{\varepsilon}\|_{L^3}\lesssim\|Y^{\varepsilon}\|_{H^{\frac{1}{3}}}\lesssim \|Y^{\varepsilon}\|_{H^1}^{\frac{2}{3}}\|Y^{\varepsilon}\|_{H^{-1}}^{\frac{1}{3}} .
$$
For any $t\leq T_{\varepsilon}$ by (\ref{4.14}) we have
\begin{equation}\label{4.17}
\begin{aligned}
\int_0^t\|Y^{\varepsilon}\|_{L^3}^3ds
&\lesssim\sup_{t\in[0,t]}\|Y^{\varepsilon}\|_{H^{-1}}\int_0^t\|Y^{\varepsilon}\|_{H^1}^2ds\\
&\lesssim\varepsilon^{\frac{\gamma_1}{2}+\gamma_2}.
\end{aligned}
\end{equation}

Then we have that for $\varepsilon$ small enough, $T_{\varepsilon}=T$, if $\gamma<\frac{\gamma_1}{2}+\gamma_2$.

Let $\gamma_1>\frac{2}{3}\gamma-\frac{2}{3}$ such that $\gamma_2=\frac{2}{3}\gamma-\frac{5}{3}$, then we only need
\begin{equation*}
\gamma_1>\frac{2}{3}\gamma+\frac{10}{3}.
\end{equation*}
i.e.
\begin{equation*}
\left\{
\begin{aligned}
&\gamma-1>\frac{2}{3}\gamma+\frac{10}{3}\\
&3(\sigma^\ast-2\delta)-1>\frac{2}{3}\gamma+\frac{10}{3}\\
&\sigma^\ast-1-2\delta+\frac{\gamma}{3}>\frac{2}{3}\gamma+\frac{10}{3}.
\end{aligned}\right.
\end{equation*}
A direct calculation yields that
\begin{equation}\label{4.18}
\left\{
\begin{aligned}
\gamma&>13\\
\sigma^\ast&>\frac{1}{3}\gamma+\frac{13}{3}+2\delta,
\end{aligned}\right.
\end{equation}
which also implies $\gamma_1=(\gamma-1)\wedge(\sigma^\ast-1-2\delta+\frac{\gamma}{3})$.

Since $R^\varepsilon=Y^{\varepsilon}+Z^{\varepsilon}$, by Lemma \ref{l4.1} we have for any $\omega\in\Omega_\delta$
\begin{equation}\label{4.19}
\begin{aligned}
\|R^\varepsilon(\omega)\|_{L^3(\mathcal{D}_T)}&\lesssim \varepsilon^{\frac{\gamma}{3}}+\varepsilon^{\sigma^\ast-2\delta}\lesssim\varepsilon^{\frac{\gamma}{3}},\\
\|R^\varepsilon(\omega)\|_{L^{\infty}(0,T;H^{-1})}^2&\lesssim\varepsilon^{\gamma-1}+\varepsilon^{\sigma^\ast-1-2\delta+\frac{\gamma}{3}}.
\end{aligned}
\end{equation}
Hence we note that
$$
w^\eps-w^\eps_A=\varepsilon\Delta(Y^\varepsilon+Z^\varepsilon)-\frac{1}{\varepsilon}\left(f(u_\varepsilon)-f(u^\varepsilon_A)\right).
$$
Using the embedding $C(\mathcal{D})\subset L^2$ we get
$$
\|\Delta(Y^\varepsilon+Z^\varepsilon)\|_{L^1(0,T;H^{-2})}\lesssim \|Y^\varepsilon\|_{L^2(0,T;L^2)}+\|Z^\varepsilon\|_{C(\mathcal{D}_T)}\lesssim\varepsilon^{\frac{\gamma}{3}}+\varepsilon^{\sigma^*-2\delta}.
$$
Similarly as above
\begin{align*}
f(u^\varepsilon)-f(u^\varepsilon_A)&=f'(u^\varepsilon_A)R^{\varepsilon}+\mathcal{N}(u^\varepsilon_A,R^\varepsilon)\\
&=f'(u^\varepsilon_A)R^\varepsilon+f''(u^\varepsilon_A+\theta(R^{\varepsilon}))(R^{\varepsilon})^2\\
&=f'(u^\varepsilon_A)R^\varepsilon+6(u^\varepsilon_A+\theta(R^{\varepsilon}))(R^{\varepsilon})^2.
\end{align*}
Since $\{u^\varepsilon_A\}$ are uniformly bounded in $\varepsilon$ and $\theta\in[0,1]$, we have that
\begin{align*}
\|f(u^\varepsilon)-f(u^\varepsilon_A)\|_{L^1(0,T;H^{-2})}
&\lesssim \|f(u^\varepsilon)-f(u^\varepsilon_A)\|_{L^1(\mathcal{D}_T)}\\
&\lesssim \|(R^\varepsilon)^3\|_{L^1(\mathcal{D}_T)}+\|(R^\varepsilon)^2\|_{L^1(\mathcal{D}_T)}+\|R^\varepsilon\|_{L^1(\mathcal{D}_T)}\\
&\lesssim \|R^\varepsilon\|_{L^3(\mathcal{D}_T)}+\|R^\varepsilon\|_{L^3(\mathcal{D}_T)}^2+\|R^\varepsilon\|_{L^3(\mathcal{D}_T)}^3\\
&\lesssim  \varepsilon^{\frac{\gamma}{3}}+\varepsilon^{\sigma^*-2\delta}\\
&\lesssim \varepsilon^{\frac{\gamma}{3}},
\end{align*}
where we use the Soblev embedding $L^1\subset H^{-2}$ in the first inequality.

Hence we deduce that
\begin{equation}\label{4.20}
\begin{aligned}
\|w^\eps-w^\eps_A\|_{L^1(0,T;H^{-2})}&\lesssim \varepsilon^{\frac{\gamma}{3}+1}+\varepsilon^{\frac{\gamma}{3}-1}\\
&\lesssim \varepsilon^{\frac{\gamma}{3}-1}.
\end{aligned}
\end{equation}
The statement of the Theorem \ref{t3.4} then follows on combining the above inequality with (\ref{4.18}).

\section{Sharp interface limit for the divergence-type noise}\label{s5}
Throughout this section we consider the singular divergence-type noise $\mathcal{W}=\nabla\cdot W$, 
where $W$ is an $L_0^2(\mathcal{D},\mathbb{R}^2)$-cylindrical Wiener process on stochastic basis $(\Omega,\mathcal{F},\mathbb{P})$. For $g\in L_0^2(\mathcal{D},\mathbb{R}^2)$, we denote its component functions by $g_1,g_2\in L_0^2(\mathcal{D})$, i.e. $g(x)=(g_1(x),g_2(x)),\forall x\in \mathcal{D}$. There exist two independent $L_0^2(\mathcal{D})$-cylindrical Wiener processes $W^1$ and $W^2$ such that $W=(W^1,W^2)$.
Similarly as in \cite{Rockner:2015uh,Rockner:2018ta}, it follows that the solution to (\ref{1.1}) 
with the divergence-type noise is distribution-valued. 
It does not appear to be possible to obtain the sharp-interface limit by directly considering (\ref{5.1}).
Thus we study the sharp interface limit of the regularized equation (\ref{1.11}) instead.

\subsection{Existence and uniqueness of solutions to equation (\ref{1.11})}

In order to consider the convolution of the noise with an approximate delta function (the standard mollifier).
we need to extend the noise to the whole space $\mathbb{R}^2$. Considering the Neumann boundary condition, it is reasonable to extend it evenly to $[-1,1]^2$ first, then do a periodical extension to the whole space. 
That is, for any function $g$ on $\mathcal{D}$ which satisfies the Neumann boundary condition, we view it as a function $\bar{g}$ on $\mathbb{R}^2$ by 
$$
\bar{g}(x):=g(|x_1+k_1|,|x_2+k_2|),\;\;\forall x=(x_1,x_2)\in\mathbb{R}^2,\;\;\forall k=(k_1,k_2)\in\mathbb{Z}^2\;\;\text{when}\;\;x+k\in[-1,1]^2.
$$
Moreover, for $x\in\mathbb{R}^2$ and $t>0$, define
$$
\bar{M}(t,x)=-\mathcal{F}^{-1}(e^{-\frac{t}{2}|\pi\cdot|^4})(x),
$$
where $\mathcal{F}^{-1}$ is the inverse Fourier transformation on $\mathbb{R}^2$. 
By Poisson summation formula, for any $(x,y)\in\mathcal{D}^2$
$$
M(t,x,y):=\sum_{k\in \mathbb{Z}^2}\left(\bar{M}(t,x+y+2k)+\bar{M}(t,x-y+2k)\right)
$$
is the kernel of $e^{-t\Delta^2}$ on $\mathcal{D}$, where $\Delta$ is the Neumann Laplacian operator on $\mathcal{D}$. A direct calculation yields that for any $g\in L^2(\mathcal{D})$
\begin{equation}\label{a.0}
\int_\mathcal{D}M(t,x,y)g(y)dy=\int_{\mathbb{R}^2}\bar{M}(t,x-y)\bar{g}(y)dy.
\end{equation}

Define
$$K(t,x,y):=-\nabla_y M(t,x,y)=\sum_{k\in \mathbb{Z}^2}\left(\bar{K}(t,x+y+2k)-\bar{K}(t,x-y+2k)\right),$$
where $\bar{K}(t,x)=(\bar{K}^1(t,x),\bar{K}^2(t,x)):=-\nabla \bar{M}(t,x)$, thus for any $t>0$, $\bar{K}^j(t,\cdot)$ is the inverse Fourier transfrmation of the function $\eta\to -\pi i\eta_je^{-\frac{t}{2}|\pi\eta|^4}$, i.e.
$$\bar{K}^j(t,x):=-\mathcal{F}^{-1}(\pi i\eta_je^{-\frac{t}{2}|\pi\eta|^4})(x).$$

We use $\mathcal{S}(\mathbb{R}^2)$ to denote the Schwartz funtion on $\mathbb{R}^2$, $\mathcal{S}'(\mathbb{R}^2)$ to denote the Schwartz distribution on $\mathbb{R}^2$ and ${}_{\mathcal{S}'(\mathbb{R}^2)}\langle\cdot,\cdot\rangle_{\mathcal{S}(\mathbb{R}^2)}$ to denote the dual between $\mathcal{S}(\mathbb{R}^2)$ and $\mathcal{S}'(\mathbb{R}^2)$.
Then we know that $\bar{K}^j(t,\cdot)\in \mathcal{S}(\mathbb{R}^2)$ for any $t>0$. Moreover we define $Z^{\varepsilon}$ by
\begin{equation}\label{a.2}
\begin{aligned}
Z^\varepsilon(t,x)&:=\varepsilon^\sigma\int_0^t\langle K(t-s,x,\cdot),dW_s\rangle_{L^2(\mathcal{D},\mathbb{R}^2)}
=\varepsilon^\sigma\sum_{j=1}^2\int_0^t{}_{\mathcal{S}'(\mathbb{R}^2)}\langle \bar{K}^j(t-s,x-\cdot),d\bar{W}^j_s\rangle_{\mathcal{S}(\mathbb{R}^2)}.
\end{aligned}
\end{equation}
Here $\bar{W}=(\bar{W}^1,\bar{W}^2)$, $\bar{W}^j,\;j=1,2$ is two i.i.d Wiener processes defined by 
$$
{}_{\mathcal{S}'(\mathbb{R}^2)}\langle\bar{W}^j,g\rangle_{\mathcal{S}(\mathbb{R}^2)}=\langle W,\tilde{g}\rangle_{L^2(\mathcal{D})},
$$
for any
$g\in\mathcal{S}(\mathbb{R}^2)$ and $\tilde{g}\in L^2(\mathcal{D})$ is defined as
$$
\tilde{g}(x):=\sum_{k\in \mathbb{Z}^2}\left(g(x+2k)-g(-x+2k)\right),\;\;x\in\mathcal{D}.
$$
For simplicity we write
$$
Z^\varepsilon(t,x)=\varepsilon^\sigma\sum_{j=1}^2\int_0^t{}_{\mathcal{S}'(\mathbb{R}^2)}\langle \bar{K}^j(t-s,x-\cdot),d\bar{W}^j_s\rangle_{\mathcal{S}(\mathbb{R}^2)}:=\varepsilon^\sigma\int_0^t{}_{\mathcal{S}'}\langle \bar{K}(t-s,x-\cdot),d\bar{W}_s\rangle_{\mathcal{S}}.
$$
We also denote
\begin{equation}
\bar{Z}^{\varepsilon}:=Z^{\varepsilon}+e^{-\varepsilon t\Delta^2}m(z),
\end{equation}
where $z\in H^{-1}$, $m(z)$ is defined in Section \ref{s2}. 
Then $\bar{Z}^\varepsilon$ is the mild solution to the linear equation
\begin{equation*}
 \left\{
   \begin{aligned}
   d\bar{Z}^\varepsilon&=-\varepsilon\Delta^2 \bar{Z}^\varepsilon+\varepsilon^{\sigma}BdW, \\
  \bar{Z}^\varepsilon(0)&\equiv m(z) \in \mathbb{R},\\
   \end{aligned}
   \right.
\end{equation*}
with Neumann boundary conditions,
\begin{equation*}
\frac{\partial \bar{Z}^\varepsilon}{\partial n}=\frac{\partial\Delta \bar{Z}^\varepsilon}{\partial n}=0\; \text{on} \;\partial \mathcal{D},
\end{equation*}
where 
  \begin{equation}
  D(B)=H^1(\mathcal{D},\mathbb{R}^2),B=\textrm{div},D(B^{\ast})=H^1(\mathcal{D}),
  B^{\ast}=-\nabla.
  \end{equation}

Let $\rho_h$ be an approximate delta function on $\mathbb{R}^2$ given by
$$\rho_h(x)=h^{-2}\rho(\frac{x}{h}), \quad\int \rho=1.$$
Define for any $(t,x)\in\mathcal{D}_T$
\begin{equation}\label{a.6}
\begin{aligned}
Z^{\varepsilon,h}(t,x):&=\varepsilon^\sigma\int_0^t{}_{\mathcal{S}'}\langle \bar{K}(\varepsilon(t-r),x-\cdot),d\bar{W}^h_s\rangle_{\mathcal{S}}\\
&=\varepsilon^\sigma\int_0^t {}_{\mathcal{S}'}\langle\bar{K}_h(\varepsilon(t-r),x-\cdot),d\bar{W}_s\rangle_{\mathcal{S}},
\end{aligned}
\end{equation}
where $\bar{W}^h=\bar{W}*\rho_h$, and $\bar{K}_h(t,x)=(\bar{K}^1_h(t,x),\bar{K}^2_h(t,x))$,
$$
\bar{K}^j_h(t,x)=\int_{\mathbb{R}^2}\bar{K}^j(t,x-y)\rho_h(y)dy.
$$

For fixed $\varepsilon,h>0$, let $\varphi^{\varepsilon,h}$ be a solution to the following equation on $\mathcal{D}$
\begin{equation}\label{a.3}
 \left\{
   \begin{aligned}
 \frac{d\varphi^{\varepsilon,h}}{dt}&=\Delta(-\varepsilon\Delta \varphi^{\varepsilon,h}+\frac{1}{\varepsilon}:f(\varphi^{\varepsilon,h}+\bar{Z}^{\varepsilon,h}):)\\
  \varphi^{\varepsilon,h}(0)&=(z-m(z))*\rho_h ,\\
   \end{aligned}
   \right.
\end{equation}
with $\Delta$ the Neumann Laplacian operator on $\mathcal{D}$. Here $:f(\varphi^{\varepsilon,h}+\bar{Z}^{\varepsilon,h}):$ is the Wick power defined by
\begin{equation}\label{a.5}
\begin{aligned}
:f(\varphi^{\varepsilon,h}+\bar{Z}^{\varepsilon,h})::=\sum_{k=0}^3C_3^k:\left(\bar{Z}^{\varepsilon,h}\right)^{3-k}:\left(\varphi^{\varepsilon,h}\right)^k
\end{aligned}
\end{equation}
where for any $k=0,1,2,3$
$$
:\left(\bar{Z}^{\varepsilon,h}\right)^{k}::=\sum_{l=0}^3C_3^k:\left({Z}^{\varepsilon,h}\right)^{k-l}:\left(e^{-\varepsilon t\Delta^2}m(z)\right)^k,
$$
\begin{align*}
:\left(Z^{\varepsilon,h}\right)^{0}:&:=1,\;:\left(Z^{\varepsilon,h}\right)::=\left(Z^{\varepsilon,h}\right),\;:\left(Z^{\varepsilon,h}\right)^{2}:=\left(Z^{\varepsilon,h}\right)^{2}-c_{h,t}^{\varepsilon}(x),\\\;:\left(Z^{\varepsilon,h}\right)^{3}:&:=\left(Z^{\varepsilon,h}\right)^{3}-3c_{h,t}^{\varepsilon}(x)\left(Z^{\varepsilon,h}\right).
\end{align*}
and
\begin{equation}\label{5.7c}
c_{h,t}^{\varepsilon}(x)=\mathbb{E}\left[Z^{\varepsilon,h}(t,x)^2\right].
\end{equation}

\begin{lemma}(\cite[Example 5.2.27]{Liu:2015vb})
For any $\varepsilon,h>0$, there exists a unique solution $\varphi^{\varepsilon,h}\in C([0,T];L^2(\mathcal{D}))$ to equation (\ref{a.3}).
\end{lemma}

Since $m(z)\in\mathbb{R}$, similar as in the proof in\cite{Mourrat:2015uo,Rockner:2015uh,Rockner:2018ta}, for any $k=1,2,3$, as $h\to 0$, $:\left(\bar{Z}^{\varepsilon,h}\right)^{k}:$ converges in $C([0,T],\mathcal{C}^\alpha)$ for any $\alpha<0$ whose limit is denoted as $:\left(\bar{Z}^{\varepsilon}\right)^{k}:$. Here $\mathcal{C}^\alpha$ is defined as the Besov space $B_{\infty,\infty}^\alpha$, see \cite{Rockner:2015uh} and the reference therein for details.
 
Then we denote
\begin{equation}\label{5.8.0}
 \left\{
   \begin{aligned}
   \frac{d\varphi^\varepsilon}{dt}&=\Delta(-\varepsilon\Delta \varphi^\varepsilon+\frac{1}{\varepsilon}:f(\varphi^\varepsilon+\bar{Z}^\varepsilon):), \\
   \varphi^\varepsilon(0)&=z-m(z) \in H_0^{-1},\\
   \end{aligned}
   \right.
\end{equation}
where
\begin{equation}\label{5.9.0}
:f(\varphi^\varepsilon+\bar{Z}^\varepsilon)::=\sum_{k=0}^3C_3^k:\left(\bar{Z}^{\varepsilon}\right)^{3-k}:\left(\varphi^{\varepsilon}\right)^k
\end{equation}

\begin{theorem}\label{t5.1}(\cite[Theorem 4.4]{Rockner:2018ta})
For $\mathbb{P}-a.s.\;\omega$,
there exists a unique solution $\varphi^\varepsilon$ to equation (\ref{5.8.0}) in $C([0,T];H_0^{-1})$ for any fixed $\varepsilon>0$.
\end{theorem}
\begin{remark}
We note that in \cite{Rockner:2018ta} the authors consider the periodical boundary condition, which is different from the Neumann boundary condition. But by our extension method as we explained before, a similar proof follows.
\end{remark}

In fact, $\varphi^\varepsilon=\lim_{h\to 0}\varphi^{\varepsilon,h}$ in $C([0,T];H_0^{-1})$.
Let $u^{\varepsilon,h}:=\varphi^{\varepsilon,h}+Z^{\varepsilon,h}$,  $u^{\varepsilon,h}$ also converges to $u^\varepsilon$ in $C([0,T];H^{-1})$, which is the unique solution to 
\begin{equation}\label{5.1}
 \left\{
   \begin{aligned}
   du^\varepsilon&=\Delta(-\varepsilon\Delta u^\varepsilon+\frac{1}{\varepsilon}:f(u^\varepsilon):)dt+\varepsilon^{\sigma}BdW, \\
   u^\varepsilon(0)&=z \in H^{-1},\\
   \end{aligned}
   \right.
\end{equation}
with Neumann boundary conditions,
\begin{equation}
\frac{\partial u^\varepsilon}{\partial n}=\frac{\partial\Delta u^\varepsilon}{\partial n}=0\; \text{on} \;\partial \mathcal{D},
\end{equation}
where $:f(u^\varepsilon):=:f(\varphi^\varepsilon+\bar{Z}^\varepsilon):$ is defined in (\ref{5.9.0}).

  \subsection{The sharp interface limit of equation (\ref{1.11})}

Similarly as in the proof of Theorem \ref{t3.4} we prove that for a suitable choice $h(\varepsilon)$, the solutions to (\ref{5.1}) will converge to the solution to deterministic Hele-Shaw model (\ref{3.2}). 

The method is a modification of the one in Section \ref{s4}. We consider the residual 
\begin{equation}\label{5.4}
R^{\varepsilon,h}:=u^{\varepsilon,h}-u^\varepsilon_A.
\end{equation}
 Let $Y^{\varepsilon,h}=R^{\varepsilon,h}-Z^{\varepsilon,h}$, which satisfies 
\begin{equation}\label{5.5}
\begin{aligned}
dY^{\varepsilon,h}=&-\varepsilon\Delta^2Y^{\varepsilon,h}dt+\frac{1}{\varepsilon}\Delta\left(f'(u^\varepsilon_A)(Y^{\varepsilon,h}+Z^{\varepsilon,h})+\mathcal{N}(u^\varepsilon_A,Y^{\varepsilon,h}+Z^{\varepsilon,h})\right)dt\\
&-\frac{c_{h,t}^{\varepsilon}}{\varepsilon}\Delta(u^\varepsilon_A+Z^{\varepsilon,h}+Y^{\varepsilon,h})+\Delta r^\varepsilon_Adt,
\end{aligned}
\end{equation}
where $c_{h,t}^{\varepsilon}$ is defined in (\ref{5.7c}).
For $Y^{\varepsilon,h}$ we also have the energy estimate:
\begin{equation}\label{5.6}
\begin{aligned}
\frac{1}{2}\frac{d\|Y^{\varepsilon,h}\|_{H^{-1}}^2}{dt}+\varepsilon\|Y^{\varepsilon,h}(t)\|_{H^1}^2=&-\frac{1}{\varepsilon}\langle f'(u^\varepsilon_A)(Y^{\varepsilon,h}+Z^{\varepsilon,h})+\mathcal{N}(u^\varepsilon_A,Y^{\varepsilon,h}+Z^{\varepsilon,h}),Y^{\varepsilon,h}\rangle\\
&-\langle r^\varepsilon_A,Y^{\varepsilon,h}\rangle+\frac{c_{h,t}^{\varepsilon}}{\varepsilon}\langle u^\varepsilon_A+Y^{\varepsilon,h}+Z^{\varepsilon,h},Y^{\varepsilon,h}\rangle.
\end{aligned}
\end{equation}

In order estimate $Y^\varepsilon$, we still need the estimation of $Z^{\varepsilon,h}$ and $c_{h,t}^{\varepsilon}$.
Analogously to Lemma \ref{l4.1} we have
\begin{lemma}\label{l5.1}
There exists a consant $C_2>0$ such that for any $0<\beta\leq 1$, 
$$
\mathbb{E}\left[\|Z^{\varepsilon,h}\|_{C(\mathcal{D}_T)}\right]\leq C_2\varepsilon^{\sigma_\ast} h^{-2},
$$
where $\sigma_\ast=\sigma-\frac{\beta}{4}$.
Then for any $\delta>0$, there exists a constant $C_\delta>0$, such that
$$\mathbb{P}\left[\Omega'_{\delta}\right]>1-C_\delta\varepsilon^\delta,$$
where $\Omega'_{\delta}=\{\|Z^{\varepsilon,h}\|_{C(\mathcal{D}_T)}\leq C_2\varepsilon^{\sigma_\ast-2\delta}h^{-2}\}$.
\end{lemma}
\proof
We follow a similar proof as in Lemma \ref{l4.1}. A factorization formula implies that 
$$
Z^{\varepsilon,h}(t,x)=\varepsilon^\sigma\frac{\sin\pi\kappa}{\pi}\int_0^t(t-s)^{\kappa-1}\langle M(\varepsilon(t-s),x-\cdot),U^{\varepsilon,h}(s)\rangle ds,
$$
where $M(t,x,y)$ is the kernel of $e^{-t\Delta^2}$
and
$$U^{\varepsilon,h}(s,x)=\int_0^t\langle (t-r)^{-\kappa}{K}_h(\varepsilon(t-r),x,\cdot),dW_s\rangle_{L^2(\mathcal{D},\mathbb{R}^2)},$$
where ${K}_h$ is defined in (\ref{a.10}). Combined with (\ref{a.6.1}), we have that
$$
|{K}_h(\varepsilon t,x,y)|\lesssim (\varepsilon t)^{-\frac{\beta}{4}}h^{-\eta}\left(|x-y|^{-\zeta}+|x+y|^{-\zeta}\right),
$$
where $\beta,\zeta,\eta\geq 0$ and $\beta+\zeta+\eta=3$.
Similarly to (\ref{2.5})-(\ref{2.9}) we have that
\begin{equation}
\begin{aligned}
\mathbb{E}\left[\left|U^{\varepsilon,h}(s,x)\right|^{2}
\right]&\lesssim \varepsilon^{-\frac{\beta}{2}}h^{-2\eta}\int_0^s\int_{\mathcal{D}}(s-r)^{-2\kappa-\frac{\beta}{2}}\left(|x+y|^{-2\zeta}+|x-y|^{-2\zeta}\right)dyds\\
&\lesssim \varepsilon^{-\frac{\beta}{2}}h^{-2\eta}s^{1-2\kappa-\frac{\beta}{2}}|x|^{2-2\zeta},
\end{aligned}
\end{equation}
where we require that 
$$
1-2\kappa-\frac{\beta}{2}>0,\quad \zeta<1
$$
Similarly to  (\ref{2.5.0}), we have that
$$\mathbb{E}\left[\|U^{\varepsilon,h}\|_{L^{2p}(\mathcal{D}_T)}\right]\lesssim \varepsilon^{\sigma-\frac{\beta}{4}}h^{-\eta}.$$

Let  $\eta=2$ and $\kappa>0$ be small enough such that $\beta<1-\zeta<2-4\kappa$, $\zeta<1$.
Similarly as in the proof of Lemma 2.7 in \cite{DaPrato:2004fs}, we have that
\begin{equation}\label{5.9}
\begin{aligned}
\mathbb{E}\left[\|Z^{\varepsilon,h}(t)\|_{\mathcal{C}(\mathcal{D}_T)}\right]&\lesssim_T\varepsilon^\sigma\mathbb{E}\left[\|U^{\varepsilon,h}\|_{L^{2p}(\mathcal{D}_T)}\right]\\
&\lesssim \varepsilon^{\sigma-\frac{\beta}{4}}h^{-2}.
\end{aligned}
\end{equation}
The statement then follows by Chebyshev's inequality.

$\hfill\Box$
\vskip.10in

For $c_{h,t}^{\varepsilon}$, we have the following estimate:
\begin{lemma}\label{l5.5}
There exists a constant $C>0$ such that for any $(t,x)\in\mathcal{D}_T$ and any $\varepsilon, h\in(0,1)$, 
$$
|c_{h,t}^{\varepsilon}(x)|\leq -C\varepsilon^{2\sigma-1}\log{h}
$$
\end{lemma}
\proof 
Following a similar argument as in (\ref{2.6}), (\ref{2.6.1}) and (\ref{2.8}), we obtain that for all $g\in(g_1,g_2)\in L^2(\mathcal{D},\mathbb{R}^2)$
\begin{align*}
\int_{\mathcal{D}}K(t,x,y)g(y)dy
&=\int_\mathcal{D}K^1(t,x,y)g_1(y)dy+\int_\mathcal{D}K^2(t,x,y)g_2(y)dy\\
&\simeq\sum_k\left(\langle g_1,|k_1|e_k\rangle+\langle g_2,|k_2|e_k\rangle\right)
e^{- t|k|^4}e_k(x).
\end{align*}
Hence 
\begin{equation}\label{a.8}
K(t,x,y)\simeq\sum_k|k|e^{-t|k|^4}e_k(x)e_k(y).
\end{equation}
where $e_k$ is defined in (\ref{2.3}). Note that $e_k(x)e_k(y)=\frac{1}{2}\left(e_k(x-y)+e_k(x+y)\right)$. Thus we obtain 
\begin{equation}
K(t,x,y)\simeq\sum_k|k|e^{-t|k|^4}(e_k(x-y)+e_k(x+y)):=P_2(t,x-y)+P_2(t,x+y).
\end{equation}
By \cite[p282, (c)]{Stein:1972jf}, we have that for any $(t,x)\in\mathcal{D}_T$,
\begin{equation}
|P_2(t,x)|\lesssim |x|^{-3}e^{-\frac{t}{|x|^4}}\lesssim \left(t^\frac{1}{3}+|x|\right)^{-3}.
\end{equation}
Thus we obtain for any $t\in[0,T]$, $x,y\in\mathcal{D}$,
$$
|K(\varepsilon t,x,y)|\lesssim\left((\varepsilon t)^{\frac{1}{4}}+|x-y|\right)^{-3}+\left((\varepsilon t)^{\frac{1}{4}}+|x+y|\right)^{-3}.
$$

We can extend the definition of $K(t,x,y)$ for $x,y\in\mathbb{R}^2$ with the same form as in (\ref{a.8}), and denote
\begin{equation}\label{a.10}
K_h(t,x,y):=\int_{\mathbb{R}^2}\rho_h(z)K(t,x,y-z)dz.
\end{equation}
Therefore (\ref{a.6}) becomes
\begin{equation*}
Z^{\varepsilon,h}(t,x)=\varepsilon^\sigma\int_0^t\langle K_h(\varepsilon(t-r),x-\cdot),dW_s\rangle_{L^2(\mathcal{D},\mathbb{R}^2)}
\end{equation*}

Then by \cite[Lemma 10.17]{Hairer:2014hd} we have that
\begin{equation}\label{a.6.1}
|{K}_h(\varepsilon t,x,y)|\lesssim\left((\varepsilon t)^{\frac{1}{4}}+|x-y|+h\right)^{-3}+\left((\varepsilon t)^{\frac{1}{4}}+|x+y|+h\right)^{-3}.
\end{equation}

Then we have that for any $(t,x)\in\mathcal{D}_T$.
\begin{equation}\label{a.7}
\begin{aligned}
|c_{h,t}^{\varepsilon}(x)|&\leq\varepsilon^{2\sigma}\int_0^t\int_{\mathcal{D}}|{K}^{\varepsilon}_h(t-r,x,y)|^2drdy\\
&\lesssim\varepsilon^{2\sigma-1}\int_0^{t\varepsilon}\int_{\mathcal{D}}\left(r^{\frac{1}{4}}+|x-y|+h\right)^{-6}drdy+\varepsilon^{2\sigma-1}\int_0^{t\varepsilon}\int_{\mathcal{D}}\left(r^{\frac{1}{4}}+|x+y|+h\right)^{-6}drdy\\
&\lesssim -\varepsilon^{2\sigma-1}\log{h}.
\end{aligned}
\end{equation}

$\hfill\Box$
\vskip.10in

The next theorem is the main result of this section.
\begin{theorem}\label{t5.4}
Let $u^{\varepsilon,h}$ be the unique solution to (\ref{5.1}) and $u^\varepsilon_A$ be defined in Theorem \ref{t3.2} with large enough $K>0$.
For some $\theta>0$ such that $\varepsilon^{\theta}\lesssim h^2$, we assume that
\begin{equation}\label{5.16}
\left\{
\begin{aligned}
\gamma&>13,\\
\sigma &>\frac{1}{3}\gamma+\frac{13}{3}+\theta.
\end{aligned}
\right.
\end{equation}
Then there exist a generic constant $C>0$ and a constant $C_\delta>0$ for all $0<\delta<\frac{\sigma}{2}-\frac{1}{6}\gamma-\frac{13}{6}-\frac{\theta}{2}$ such that the following estimates hold
\begin{equation}\label{5.17}
\begin{aligned}
\mathbb{P}\left[\|R^{\varepsilon,h}\|_{L^3(\mathcal{D}_T)}\leq C\varepsilon^{\frac{\gamma}{3}}\right]\geq 1-C_\delta\varepsilon^\delta,\\
\mathbb{P}\left[\|R^{\varepsilon,h}\|_{L^{\infty}(0,T;H^{-1})}^2\leq C\left(\varepsilon^{\gamma-1}+\varepsilon^{\sigma_\ast-1-2\delta-\theta}\right)\right]\geq 1-C_\delta\varepsilon^\delta,\\
\mathbb{P}\left[\|v^{\varepsilon,h}-w^\eps_A\|_{L^{1}(0,T;H^{-2})}^2\leq C\varepsilon^{\frac{\gamma}{3}-1}\right]\geq 1-C_\delta\varepsilon^\delta.
\end{aligned}
\end{equation}
\end{theorem}

\proof 
We proceed similarly as in Section~\ref{s4}.
We define a stopping time
\begin{equation}\label{5.10}
T^{\varepsilon,h}:=T\wedge\inf\{t>0:\int_0^t\|Y^{\varepsilon,h}(\tau)\|_{L^3}^3d\tau>\varepsilon^{\gamma}\}. 
\end{equation}
Then let $t<T^{\varepsilon,h}$ and fix an $\omega\in \Omega'_\delta$. Since
\begin{equation}\label{5.11}
h^{-2}\lesssim\varepsilon^{-\theta}
\end{equation}
for some $\theta>0$. We have that
\begin{equation}\label{5.12}
-\log h\lesssim -\frac{\theta}{2}\log \varepsilon\lesssim\varepsilon^{-\delta},\quad |c_{h,t}^{\varepsilon}|\lesssim\varepsilon^{2\sigma-1-\delta}.
\end{equation}

For $\frac{c_{h,t}^{\varepsilon}}{\varepsilon}\langle u^\varepsilon_A+Y^{\varepsilon,h}+Z^{\varepsilon,h},Y^{\varepsilon,h}\rangle$ we have that for small enough $\varepsilon$
\begin{equation*}
\int_0^t\frac{c_{h,t}^{\varepsilon}}{\varepsilon}|\langle u^\varepsilon_A+Y^{\varepsilon,h}+Z^{\varepsilon,h},Y^{\varepsilon,h}\rangle |d\tau\lesssim \varepsilon^{\frac{\gamma}{3}-1}c_{h,t}^{\varepsilon}\lesssim \varepsilon^{2\sigma+\frac{\gamma}{3}-2-\delta}.
\end{equation*}

For the rest terms on the right hand side of (\ref{5.6}), we follow the proof in Section \ref{s4} by repalcing the estimate for $Z^{\varepsilon}$ with the estimate of $Z^{\varepsilon,h}$ in Lemma \ref{l5.1}. Thus we have that for small enough $\varepsilon$ and $t\leq T^{\varepsilon,h}$
\begin{equation*}
\sup_{\tau\in[0,t]}\|Y^{\varepsilon,h}(\tau)\|_{H^{-1}}^2d\tau\lesssim \varepsilon^{2\sigma+\frac{\gamma}{3}-2-\delta}+\varepsilon^{\gamma-1}+\varepsilon^{\sigma_\ast-1-2\delta+\frac{\gamma}{3}-\theta}+\varepsilon^{3(\sigma_\ast-2\delta-\theta)-1}.
\end{equation*}
Also,
\begin{equation*}
\int_0^t\|Y^{\varepsilon,h}(\tau)\|_{H^1}^2d\tau\lesssim \varepsilon^{2\sigma+\frac{\gamma}{3}-3-\delta}+\varepsilon^{\frac{2}{3}(\gamma-1)-1}+\varepsilon^{\sigma_\ast-2-2\delta+\frac{\gamma}{3}-\theta}+\varepsilon^{3(\sigma_\ast-2\delta-\theta)-2}+\varepsilon^{\gamma-2}.
\end{equation*}

Hence we have
\begin{equation}
\sup_{\tau\in[0,t]}\|Y^{\varepsilon,h}(\tau)\|_{H^{-1}}^2d\tau\lesssim\varepsilon^{\gamma_1},\quad \int_0^t\|Y^{\varepsilon,h}(\tau)\|_{H^1}^2d\tau\lesssim\varepsilon^{\gamma_2},
\end{equation}
where 
$$\gamma_1:=(2\sigma+\frac{\gamma}{3}-2-\delta)\wedge(\sigma_\ast-1-2\delta-\theta+\frac{\gamma}{3})\wedge(3(\sigma_\ast-2\delta-\theta)-1)\wedge(\gamma-1),$$
$$\gamma_2:=(\gamma_1-1)\wedge(\frac{2}{3}\gamma-\frac{5}{3}).$$
Similarly to (\ref{4.17}), we have
\begin{equation*}
\int_0^t\|Y^{\varepsilon,h}\|_{L^3}^3d\tau\lesssim\varepsilon^{\frac{\gamma_1}{2}+\gamma_2}.
\end{equation*}
In order to prove $T^{\varepsilon,N}=T$ for small enough $\varepsilon$, we need to prove $\gamma<\frac{1}{2}\gamma_1+\gamma_2$.  
First we assume that $\gamma_2=\frac{2}{3}\gamma-\frac{5}{3}$, i.e. 
\begin{equation*}
\gamma_1>\frac{2}{3}\gamma-\frac{2}{3}.
\end{equation*}
Then $\frac{1}{2}\gamma_1+\gamma_2>\gamma$ yields
\begin{equation*}
\gamma_1>\frac{2}{3}\gamma+\frac{10}{3}.
\end{equation*}
A direct calculation yields that
\begin{equation*}
\left\{
\begin{aligned}
\gamma&>13,\\
\sigma &>\frac{1}{3}\gamma+\frac{13}{3}+2\delta+\theta+\frac{\beta}{4},
\end{aligned}
\right.
\end{equation*}
which implies that
$$
\gamma_1=(\sigma_\ast-1-2\delta-\theta+\frac{\gamma}{3})\wedge(\gamma-1).
$$
Since $\delta,\beta>0$ can be arbitrarily small, we can only assume that (\ref{5.16}) hold and let $0<2\delta<\sigma-\frac{1}{3}\gamma-\frac{13}{3}-\theta$.

Since $R^{\varepsilon,h}=Y^{\varepsilon,h}+Z^{\varepsilon,h}$, and $H^1\subset L^3$, we can obtain the estimate of $R^{\varepsilon,h}$ which is similar to (\ref{4.19}). Moreover let
$$
v^{\varepsilon,h}:=-\varepsilon\Delta u^{\varepsilon,h}+\frac{1}{\varepsilon}\left(f(u^{\varepsilon,h})-3c_{h,t}^{\varepsilon}u^{\varepsilon,h}\right),
$$
similarly to (\ref{4.20}), we obtain that
\begin{align*}
\|v^{\varepsilon,h}-w^\eps_A\|_{L^{1}(0,T;H^{-2})}^2&
\lesssim\varepsilon\|R^{\varepsilon,h}\|_{L^2(\mathcal{D}_T)}+\frac{1}{\varepsilon}\|R^{\varepsilon,h}\|_{L^1(\mathcal{D}_T)}+\frac{1}{\varepsilon}\|c_{h,t}^{\varepsilon}u^{\varepsilon,h}\|_{L^{1}(0,T;H^{-2})}^2\\
&\lesssim \varepsilon\|R^{\varepsilon,h}\|_{L^3(\mathcal{D}_T)}+\frac{1}{\varepsilon}\|R^{\varepsilon,h}\|_{L^3(\mathcal{D}_T)}+\frac{1}{\varepsilon}\|c_{h,t}^{\varepsilon}\|_{L^\infty}\|R^{\varepsilon,h}+u^\varepsilon_A\|_{L^3(\mathcal{D}_T)}\\
&\lesssim\varepsilon^{\frac{\gamma}{3}-1}-\varepsilon^{2\sigma-2}\log h\lesssim\varepsilon^{\frac{\gamma}{3}-1}.
\end{align*}

$\hfill\Box$
\vskip.10in

\begin{remark}
It is easy to see that (\ref{5.16}) implies the condition $\sigma>\frac{26}{3}+\theta$. 
This implies that the slower $h$ converges to $0$ than $\varepsilon$, the smaller $\sigma$ could be.
Since $\theta$ can be arbitrarily small, the lower bound for $\sigma$ is $\frac{26}{3}$.
\end{remark}

The next corollary 
is a simple consequence of Theorem~\ref{t5.4} and can be shown as Corollary~\ref{c3.5}.
\begin{corollary}
There exist a subsequence $\{\varepsilon_k\}_{k=1}^\infty$ and $\{h_k\}_{k=1}^\infty$ with $\varepsilon_k^\theta\lesssim h^2_k$ such that for $\mathbb{P}-a.s.\;\omega\in\Omega$
$$
\lim_{k\to\infty}u^{\varepsilon_k,h_k}={ 1-2\chi_{\mathcal{D}_t^-}}\;\;\text{in}\;\;L^3({\mathcal{D}_{T}}),
$$
where $\mathcal{D}_t^-$ is the interior of $\Gamma_t$ in $\mathcal{D}$.
\end{corollary}
%\proof
%The proof is the same as Corollay \ref{c3.5}, we ignore it here for simplicity.

 \appendix
  \renewcommand{\appendixname}{Appendix~\Alph{section}}

\bibliographystyle{alpha}%

\end{document}